\documentclass[10pt]{article}
\usepackage{amsfonts,amsthm, amsmath,cite}
\usepackage{mathtools}
\usepackage{epsfig}
\usepackage{makeidx}
\usepackage{graphicx,epstopdf}
\DeclareGraphicsExtensions{.ps,.png,.jpg}
\usepackage{enumerate}

\makeindex
\usepackage{subfig}
\usepackage{color}
\usepackage{bm}
\usepackage{setspace}
\def\E{{\mathrm{E}}}

\let\originalleft\left
\let\originalright\right
\def\left#1{\mathopen{}\originalleft#1}
\def\right#1{\originalright#1\mathclose{}}

\newcommand{\ncom}{\newcommand}
\ncom{\ul}{\underline}
\ncom{\beq}{\begin{equation}}
\ncom{\eeq}{\end{equation}}
\ncom{\bea}{\begin{eqnarray*}}
\ncom{\eea}{\end{eqnarray*}}
\ncom{\beqa}{\begin{eqnarray}}
\ncom{\eeqa}{\end{eqnarray}}
\ncom{\nno}{\nonumber}
\ncom{\non}{\nonumber}
\ncom{\ds}{\displaystyle}
\ncom{\half}{\frac{1}{2}}
\ncom{\mbx}{\makebox{.25cm}}
\ncom{\hs}{\mbox{\hspace{.25cm}}}
\ncom{\rar}{\rightarrow}
\ncom{\Rar}{\Rightarrow}
\ncom{\noin}{\noindent}
\ncom{\bc}{\begin{center}}
\ncom{\ec}{\end{center}}
\ncom{\sz}{\scriptsize}
\ncom{\rf}{\ref}
\ncom{\s}{\sqrt{2}}
\ncom{\sgm}{\sigma}
\ncom{\Sgm}{\Sigma}
\ncom{\psgm}{\sigma^{\prime}}
\ncom{\dt}{\delta}
\ncom{\Dt}{\Delta}
\ncom{\lmd}{\lambda}
\ncom{\Lmd}{\Lambda}
\ncom{\Th}{\Theta}
\ncom{\e}{\eta}
\ncom{\eps}{\epsilon}
\ncom{\pcc}{\stackrel{P}{>}}
\ncom{\lp}{\stackrel{L_{p}}{>}}
\ncom{\dist}{{\rm\,dist}}
\ncom{\sspan}{{\rm\,span}}
\ncom{\re}{{\rm Re\,}}
\ncom{\im}{{\rm Im\,}}
\ncom{\sgn}{{\rm sgn\,}}
\ncom{\ba}{\begin{array}}
\ncom{\ea}{\end{array}}
\ncom{\hone}{\mbox{\hspace{1em}}}
\ncom{\htwo}{\mbox{\hspace{2em}}}
\ncom{\hthree}{\mbox{\hspace{3em}}}
\ncom{\hfour}{\mbox{\hspace{4em}}}
\ncom{\vone}{\vskip 2ex}
\ncom{\vtwo}{\vskip 4ex}
\ncom{\vonee}{\vskip 1.5ex}
\ncom{\vthree}{\vskip 6ex}
\ncom{\vfour}{\vspace*{8ex}}
\ncom{\norm}{\|\;\;\|}
\ncom{\integ}[4]{\int_{#1}^{#2}\,{#3}\,d{#4}}
\ncom{\vspan}[1]{{{\rm\,span}\{ #1 \}}}
\ncom{\dm}[1]{ {\displaystyle{#1} } }
\ncom{\ri}[1]{{#1} \index{#1}}

\newtheorem{theorem}{\bf Theorem}[section]
\newtheorem{remark}{\bf Remark}[section]

\newtheorem{lemma}{Lemma}[section]

\newtheorem{example}{Example}[section]
\newtheorem{definition}{Definition}[section]
\newtheoremstyle
    {remarkstyle}
    {}
    {11pt}
    {}
    {}
    {\bfseries}
    {:}
    {     }
    {\thmname{#1} \thmnumber{#2} }

\theoremstyle{remarkstyle}

\def \R{{{\rm I{\!}\rm R}}}

\usepackage[round]{natbib}
\begin{document}

\newpage

\begin{center}
{\Large \bf Cross-codifference for bidimensional VAR(1) models\\ with infinite variance}
\end{center}
\vone
\begin{center}
{\large Aleksandra Grzesiek$^1$, Marek Teuerle and Agnieszka Wy{\l}oma{\'n}ska} \\
{{\emph{Faculty of Pure and Applied Mathematics},\\ \emph{Hugo Steinhaus Center,\\ Wroc{\l}aw University of Science and  Technology,}\\
			 \emph{Wybrze\.ze Wyspia{\'n}skiego 27, 50-370 Wroc{\l}aw, Poland}}}
\end{center}
\begin{center}
	$^1$\emph{Corresponding author:} aleksandra.grzesiek@pwr.edu.pl
\end{center}
\vtwo

\noindent{\bf Abstract}:
In this paper we consider the problem of a measure that allows us to describe the spatial and temporal dependence structure of multivariate time series with innovations having infinite variance. By using recent results obtained in the problem of temporal dependence structure of univariate stochastic processes, where the auto-codifference was used, we extend its idea and propose a cross-codifference measure for a general vector autoregressive model of order 1 (VAR(1)).
Next, we derive an analytical results for VAR(1) model with Gaussian and sub-Gaussian innovations, that are characterized by finite and infinite variance, respectively. We emphasize that obtained expressions perfectly agree with the empirical counterparts. Moreover, we show that for the considered processes the cross-codifference  simplifies to the well-established cross-covariance measure in case of Gaussian white noise. Last part of the work is devoted to the statistical estimation of VAR(1) parameters based on the empirical cross-codifference. Again, we demonstrate via Monte Carlo simulations that proposed methodology works correctly.
\vone \noindent{\it Key words:} dependence measure, cross-covariance, codifference, infinite variance, sub-Gaussian
\vone \noindent{\it Short version of title:} Cross-codifference for bidimensional VAR(1) models
\vtwo
\setcounter{equation}{0}
\section{Introduction}
In statistics, the cross-covariance (called called cross-correlation) is a measure of similarity between two stochastic processes, commonly used to find features in an unknown process by comparing it to a known one. It is a function of the relative time between examined systems. The cross-covariance statistics is extremely useful in structure of dependence investigation in case of bidimensional process. It gives information about the relationship between components of the bidimensional process and takes under consideration their delay in time. The cross-covariance may be also useful in the problem in model identification \citep{model}. It is especially popular in the time series analysis and signal processing \citep{c1,c2,c3,c4,c5,c6}.\\
\noindent However, for processes with infinite variance the theoretical cross-covariance does not exist, therefore, it can not be considered as the appropriate measure of similarity between components for bidimensional processes. In this case the alternative measures have to be used. In one-dimensional case one can find in the literature different alternative measures of dependence that can replace the classical autocovariance in infinite variance case. Here, we mention only the covariation \citep{cova}, fractional order covariance \citep{floc1,floc2,floc3,floc4} or codifference \citep{Taqqu,nowicka1}. Especially the last measure found many interesting applications \citep{Taqqu,cod1,cod2,Chechkin2015,nowicka2}. It is defined for any infinitely divisible processes and its estimator have relatively simple form \citep{cod1}. 
The other mentioned alternative measures are defined only for special processes, namely $\alpha-$stable-based, therefore there is a limited number of their applications.\\
\noindent In this paper we fill the gap in the description of the dependence structure for bidimensional processes with infinite variance and introduce the new measure called cross-codifference. We show it can be considered as the extension of the cross-covariance since it reduces to the classical measure in the Gaussian case. We consider here the bidimensional vector autoregressive model of order 1 (VAR(1)) where the noise is modeled by sub-Gaussian process belonging to the class of $\alpha-$stable distributions. The sub-Gaussian processes are well-known models and one can find their various applications \citep{sub1,sub2,sub3}, for example in finance to model stock indices returns \citep{sub4}. We consider also the bidimensional VAR(1) model with Gaussian noise in order to show the introduced cross-codifference is universal, it can be applied to finite and infinite variance processes and gives the same message about the similarity between components of bidimensional process and their delay in time. \\
\noindent For both VAR(1) models we calculate the cross-codifference in general case and pay attention on exemplary values of their parameters. Finally we introduce the estimator of the new measure and check its effectiveness for simulated trajectories. As the possible application of the obtained theoretical results we demonstrate the new estimation method for VAR(1) model parameters.
\section{General bidimensional VAR(1) model and the cross-covariance measure}
In this section we introduce the general VAR(1) model and  give its important characteristics used in the further analysis.
\begin{definition}
	\textup{\citep{Brockwell}} The time series $\left\{\mathbf{X}\left(t\right)\right\}=\left\{X_1\left(t\right),X_2\left(t\right)\right\}$ is a bidimensional VAR(1) model if $\left\{\mathbf{X}\left(t\right)\right\}$ is weak-sense stationary and if for every $t$ it satisfies the following equation
	\begin{align}
	\mathbf{X}\left(t\right)-\Theta\ \mathbf{X}\left(t-1\right)=\mathbf{Z}\left(t\right),
	\label{eq1}
	\end{align}
	where $\{\mathbf{Z}\left(t\right)\}$ is a bidimensional white noise. 
\end{definition} 
\noindent We remind, the $\left\{\mathbf{Z}\left(t\right)\right\}$ is a bidimensional white noise if it is weak-sense stationary with mean vector $\mathbf{0}$ and the covariance matrix function
\begin{align}
\Gamma\left(h\right)&=\left\{ \begin{array}{ll}\Sigma & \mbox{if}~h=0,\\
0& \mbox{otherwise}.
\end{array} \right.
\label{eq2}
\end{align}
\noindent If we use the vector notation
\begin{align*}
\mathbf{X}\left(t\right)=\left[ \begin{array}{ccc}X_1\left(t\right)\\X_2\left(t\right)\end{array}\right],
\qquad
\Theta= \left[
\begin{array}{cc}
a_1 & a_2\\
a_3 & a_4
\end{array}
\right],
\qquad
\mathbf{Z}(t)=\left[ \begin{array}{ccc}Z_1\left(t\right)\\Z_2\left(t\right)\end{array}\right],
\end{align*}
we can rewrite (\ref{eq1}) as
\begin{align*}
\left[ \begin{array}{ccc}X_1\left(t\right)\\X_2\left(t\right)\end{array}\right] - \left[
\begin{array}{cc}
a_1 & a_2\\
a_3 & a_4
\end{array}
\right] \left[ \begin{array}{ccc}X_1\left(t-1\right)\\X_2\left(t-1\right)\end{array}\right] =
\left[ \begin{array}{ccc}Z_1\left(t\right)\\Z_2\left(t\right)\end{array}\right],
\end{align*}
what is equivalent to the following system of recursive equations:
\begin{align}
\left\{ \begin{array}{ll}
X_1\left(t\right)-a_1\ X_1\left(t-1\right)-a_2\ X_2\left(t-1\right)= Z_1\left(t\right),\\
X_2\left(t\right)-a_3\ X_1\left(t-1\right)-a_4\ X_2\left(t-1\right)= Z_2\left(t\right).
\end{array} \right.
\label{eq2}
\end{align}
Moreover, for a bidimensional VAR(1) process under certain condition, that will be clarified later, we can express $\mathbf{X}\left(t\right)$ given in (\ref{eq1}) as, \citep{Brockwell}
\begin{align}
\mathbf{X}(t)=\sum_{j=0}^{\infty}\Theta^j\mathbf{Z}\left(t-j\right).
\label{eq3}
\end{align}
Assuming that
\begin{align*}
\Theta^j= \left[
\begin{array}{cc}
a_1^{\left(j\right)} & a_2^{\left(j\right)}\\
a_3^{\left(j\right)} & a_4^{\left(j\right)}
\end{array}
\right],
\end{align*}
one can rewrite Eq. (\ref{eq3}) as
\begin{align*}
\left[ \begin{array}{ccc}X_1\left(t\right)\\X_2\left(t\right)\end{array}\right]= \sum_{j=0}^{\infty}\left[
\begin{array}{cc}
a_1^{\left(j\right)} & a_2^{\left(j\right)}\\
a_3^{\left(j\right)} & a_4^{\left(j\right)}
\end{array}
\right] \left[ \begin{array}{ccc}Z_1\left(t-j\right)\\Z_2\left(t-j\right)\end{array}\right],
\end{align*}
what is equivalent to
\begin{align}
\left\{ \begin{array}{ll}
X_1\left(t\right) = \sum_{j=0}^{\infty} \left(a_1^{\left(j\right)} Z_1\left(t-j\right) + a_2^{\left(j\right)} Z_2\left(t-j\right)\right),\vspace{1ex}\\
X_2\left(t\right) = \sum_{j=0}^{\infty} \left(a_3^{\left(j\right)} Z_1\left(t-j\right) + a_4^{\left(j\right)} Z_2\left(t-j\right)\right).
\end{array} \right.
\label{eq4}
\end{align}
In this paper, in the general case we consider the bidimensional VAR(1) model with bidimensional sub-Gaussian distribution. Moreover, we also concentrate on the Gaussian case, namely when the white noise $\{\mathbf{Z}\left(t\right)\}$ in equation (\ref{eq1}) is Gaussian. We should mention, the equation (\ref{eq4}) is satisfied if the $\{\mathbf{X}\left(t\right)\}$ model is bounded in the sense of the norm for appropriate space of random variables. Therefore there is a need to prove the conditions that guarantee the boundary solution of equation (\ref{eq1}). In the considered cases, namely Gaussian and sub-Gaussian, the norms are different therefore the boundary conditions will be proved separately for each case, see sections \ref{var_gauss} and \ref{var_subgauss}.

\subsection{Bidimensional Gaussian noise}\label{var_gauss}
In the first case, we consider bidimensional VAR(1) model with Gaussian noise. Let us denote a random vector in $\mathbb{R}^d$  as 
\begin{equation}
\mathbf{G}=\left(G_1,G_2,\ldots,G_d\right).
\label{gvector}
\end{equation}
The vector $\mathbf{G}$ is called to be a $d$-dimensional Gaussian vector if every linear combination of its components $Y= a_1G_1 + \ldots + a_dG_d$ has Gaussian distribution. That is, for any constant vector $\mathbf{a} \in \mathbb{R}^d$ the random variable $Y = \mathbf{a}^{T}X$ has a (one-dimensional) Gaussian distribution. The characteristic function of Gaussian random vector $\mathbf{G}$ takes the form \citep{feller}
\begin{equation}
\phi_{\mathbf{G}}\left(\bm{\theta}\right)=\phi_{\mathbf{G}}\left(\left(\theta_1,\theta_2,\ldots,\theta_d\right)\right)=\exp\left\{i\bm{\theta}^{T}\bm{\mu}-\frac{1}{2}\bm{\theta}^{T}\Sigma\bm{\theta}\right\},
\label{characteristic2}
\end{equation}
where $\bm{\mu}$ is $d$-dimensional mean vector $\bm{\mu}=\left(\E G_1,\E G_2,\ldots,\E G_d\right)^{T}$ and $\Sigma$ is $d\times d$ covariance matrix $\Sigma=\left[\mathrm{Cov}\left(G_i,G_j\right);i,j=1,\ldots,d\right]$. This extends to the Gaussian process as follows. A stochastic process $\left\{G(t), t \in T\right\}$ is Gaussian if and only if its finite dimensional sets, $\left(G\left(t_1\right),G\left(t_2\right),\ldots,G\left(t_d\right)\right)$, $d \ge 1$, are Gaussian random vectors, i.e. every finite linear combination of them has Gaussian distribution.\\

\noindent In this paper, for the Gaussian VAR(1) model we take bidimensional Gaussian white noise given by (\ref{gvector}):
\begin{equation}
\mathbf{Z}\left(t\right)=\left[ \begin{array}{ccc}Z_1\left(t\right)\\Z_2\left(t\right)\end{array}\right]=\left[ \begin{array}{ccc}G_1\\G_2\end{array}\right],
\label{eq7}
\end{equation}
where $G=\left(G_1,G_2\right)$ is a zero-mean Gaussian vector in $\mathbb{R}^2$. Therefore, from the formula (\ref{characteristic2}), the characteristic function of $\mathbf{Z}\left(t\right)$ has the following form
\begin{align}
\phi_{\mathbf{Z}(t)}\left(\theta\right)=\E\left(\exp\left\{i\theta_1 Z_1\left(t\right)+i \theta_2Z_2\left(t\right)\right\}\right)=\exp\left\{-\frac{1}{2}\left(\theta_1^2 R_{11}+2\theta_1\theta_2R_{12}+\theta_2^2R_{22}\right)\right\},
\label{eq8}
\end{align}
where $R_{11}=\E G_1^2$, $R_{12}=\E\left(G_1G_2\right)$ and $R_{22}=\E G_2^2$.\\
\begin{remark} The bidimensional VAR(1) model with Gaussian noise defined in (\ref{eq1}) has bounded solution given by formula (\ref{eq4}) if the following conditions hold 
\begin{eqnarray}\label{r1}
\sum_{j=0}^{\infty}\left(|a_1^{(j)}|+ |a_2^{(j)}|\right)<\infty\quad\mathrm{and}\quad
\sum_{j=0}^{\infty} \left(|a_3^{(j)}|+ |a_4^{(j)}|\right)<\infty.
\end{eqnarray}
\end{remark}
\noindent Proof: The proof follows directly from formula (\ref{eq4}), and the properties of the $L^2$ norm for second order random variables,   namely
\[\|X_1(t)\|_2={\left\|\sum_{j=0}^{\infty} \left(a_1^{(j)} Z_1(t-j) + a_2^{(j)} Z_2(t-j)\right)\right\|}_{2}\leq \sum_{j=0}^{\infty} \left(\left|a_1^{(j)}\right| R_{11}+ \left|a_2^{(j)}\right| R_{22}\right).\]
The above is finite if $\sum_{j=0}^{\infty} \left(\left|a_1^{(j)}\right|+\left|a_2^{(j)}\right|\right)<\infty$. By the same reasoning we obtain second condition that is associated with the process $\{X_2(t)\}$. Thus we obtain the statement.
\begin{flushright}
$\Box$
\end{flushright}
In the further analysis we assume the parameters satisfy conditions (\ref{r1}).
\subsection{Bidimensional sub-Gaussian noise}\label{var_subgauss}
In the second case, we consider bidimensional VAR(1) model with  sub-Gaussian noise. Because in the definition of sub-Gaussian distribution and process there appears a notion of $\alpha-$stable random variable, thus we remind one of the definitions of $\alpha-$stable random variable. Important properties and facts related to class of $\alpha-$stable random variables one can find the literature, see for instance \cite{Taqqu}.
\begin{definition}We say a random variable $A$ has $\alpha-$stable distribution with parameters $\alpha,\sigma,\beta$ and $\mu$ ($S_{\alpha}\left(\sigma,\beta,\mu\right)$) if its characteristic function is given by
\begin{eqnarray}\label{data}
 \phi_A(\theta)=
\left\{
\begin{array}{ll} 
\exp\left\{-\sigma^{\alpha}|\theta|^{\alpha}\left\{1-i\beta \mathrm{sign}(\theta)\tan\left(\pi\alpha/2\right)\right\}+i\mu \theta\right\}& \mbox{for $\alpha\neq 1$,}\\
&\\
\exp\left\{-\sigma|\theta|\{1+i\beta \mathrm{sign}(\theta)\frac{2}{\pi}\log(|\theta|\}+i\mu \theta\right\}& \mbox{for $\alpha=1$.}
\end{array}
\right.
\end{eqnarray}
\end{definition}
\noindent In the above definition the $\alpha\in (0,2]$ parameter is called the stability index, $\sigma>0$ - scale parameter, $\beta\in [-1,1]$ - skewness and $\mu\in \R$ is the shift parameter. We say that random variable has symmetric $\alpha-$stable distribution (around zero) if $\beta=\mu=0$. We refer the readers to classical literature of $\alpha-$stable distributed random variables and processes \citep{Taqqu,nolan}. The $\alpha-$stable distribution has many interesting properties and therefore has found various applications. However, one of the disadvantage is the infinite variance for most of the cases (except $\alpha=2$, Gaussian case). It raises many problems for instance in the estimation, statistical investigation and description of the dependence structure of $\alpha-$stable-based models. This problem is also highlighted in the current paper. After this remark we can define the sub-Gaussian random variables.

\noindent Let us consider a zero-mean Gaussian random variable $G$ and an $\alpha/2-$stable totally skewed to the right (i.e. for $\beta=1$) random variable $A$, and assume $G$ and $A$ to be independent. Then a random variable  $$Z=A^{1/2}G$$ has symmetric $\alpha-$stable distribution \citep{Taqqu}. This result extends to random vectors as follows. If we take a random variable 
$$A \sim S_{\alpha/2}\left(\sigma=\cos{\left(\frac{\pi\alpha}{4}\right)}^{2/\alpha},1,0\right)$$ 
and a zero-mean Gaussian vector in $\mathbb{R}^d$ $$\mathbf{G}=\left(G_1,G_2,\ldots,G_d\right),$$ and $G$ is independent from $A$, then the random vector
\begin{equation}
\mathbf{Z}=\left(A^{1/2}G_1,A^{1/2}G_2,\ldots,A^{1/2}G_d\right)
\label{subGvector}
\end{equation}
has symmetric $\alpha-$stable distribution in $\mathbb{R}^d$ and is called a $d$-dimensional sub-Gaussian vector with underlying Gaussian vector $\mathbf{G}$ \citep{Taqqu}. The characteristic function of sub-Gaussian random vector takes the following form
\begin{equation}
\phi_\mathbf{Z}\left(\theta\right)=\phi_\mathbf{Z}\left(\left(\theta_1,\theta_2,\ldots,\theta_d\right)\right)=\exp\left\{-{\left|\frac{1}{2}\sum_{i=1}^{d}\sum_{j=1}^{d}\theta_i\theta_jR_{ij}\right|}^{\alpha/2}\right\},
\label{characteristic}
\end{equation}
where $R_{ij}=\E\left(G_iG_j\right)$ is a covariance between $G_i$ and $G_j$ \citep{Taqqu}. 
We can extend the definition of bidimensional sub-Gaussian distribution  to the sub-Gaussian process.

\noindent If we take a Gaussian process $\left\{G\left(t\right), t \in T\right\}$ and an $\alpha/2-$stable totally skewed to the right random variable $A$, and assume that $A$ and $\left\{G\left(t\right), t \in T\right\}$ are independent, then the process $\left\{Z\left(t\right)=A^{1/2}G\left(t\right), t \in T\right\}$ is called a sub-Gaussian process with underlying Gaussian process $\left\{G\left(t\right), t \in T\right\}$. Its finite dimensional sets, $\left(Z\left(t_1\right),Z\left(t_2\right),\ldots,Z\left(t_d\right)\right)$, $d \ge 1$, constitute the sub-Gaussian random vector introduced in (\ref{subGvector}) \citep{Taqqu}. 

In this paper, for the sub-Gaussian VAR(1) model we take a two-dimensional sub-Gaussian vector introduced in (\ref{subGvector}):
\begin{equation}
\mathbf{Z}(t)=\left[ \begin{array}{ccc}Z_1(t)\\Z_2(t)\end{array}\right]=\left[ \begin{array}{ccc}A^{1/2}G_1\\A^{1/2}G_2\end{array}\right],
\label{eq10}
\end{equation}
where $G=(G_1,G_2)$ is a zero-mean Gaussian vector in $\mathbb{R}^2$ and $A \sim S_{\alpha/2}\left(\sigma=\cos\left(\frac{\pi\alpha}{4}\right)^{2/\alpha},1,0\right)$. Moreover, \mbox{$G$ and $A$} are independent. From the formula (\ref{characteristic}), the characteristic function of $\mathbf{Z}(t)$ has the form
\begin{equation}
\phi_{\mathbf{Z}(t)}(\theta)=\E\left(\exp\{i\theta_1 Z_1(t)+i \theta_2Z_2(t)\}\right)=\exp\left\{-{\left(\frac{1}{2}\right)}^{\alpha/2}\left|\theta_1^2 R_{11}+2\theta_1\theta_2R_{12}+\theta_2^2R_{22}\right|^{\alpha/2}\right\},
\label{eq11}
\end{equation}
where $R_{11}=\E G_1^2$, $R_{12}=\E(G_1G_2)$ and $R_{22}=\E G_2^2$. In this case $\{\mathbf{Z}(t)\}$ is called a zero-mean sub-Gaussian white noise.

Since for each $t$ the components of  $\mathbf{Z}(t)$, namely $Z_1(t)$ and $Z_2(t)$, have\, symmetric $\alpha-$stable distribution, therefore in order to prove the conditions that guarantee existence of bounded solution (given by (\ref{eq4})) of equation (\ref{eq1}) we need first to introduce a norm in the space of symmetric $\alpha-$stable random variables. Here we consider only case $\alpha>1$.
\begin{definition}\textup{\citep{Taqqu}} Let $A_1$ and $A_2$ be two symmetric $\alpha-$stable random variables and $\alpha>1$, then the covariation between $A_1$ and $A_2$ is defined as
\begin{eqnarray}
[A_1,A_2]_{\alpha}=\int_{S_1}s_1s_2^{\langle\alpha-1\rangle}\Gamma(ds),
\end{eqnarray}
where $|a|^{\langle p\rangle}=|a|^p \mathrm{sign}(a)$, $\Gamma(\cdot)$ is a spectral measure of vector $(A_1,A_2)$ and $S_1$ is a unit circle.
\end{definition}
\begin{definition}\textup{\citep{Taqqu}} In the space of symmetric $\alpha-$stable random variables with $\alpha>1$, a covariation norm for a random variable $A$ from this space is defined as
\begin{eqnarray}
\|A\|_{\alpha}=([A,A]_{\alpha})^{1/\alpha}.
\end{eqnarray} 
\end{definition} 
\begin{remark}
The bidimensional VAR(1) model defined in (\ref{eq1}) with sub-Gaussian noise  has bounded solution given by formula (\ref{eq4}) if the following conditions hold 
\begin{equation}\label{r2}
\sum_{j=0}^{\infty}\left(\left|a_1^{(j)}\right|+\left|a_2^{(j)}\right|\right)<\infty\quad \mathrm{and}\quad
\sum_{j=0}^{\infty}\left(\left|a_3^{(j)}\right|+\left|a_4^{(j)}\right|\right)<\infty.
\end{equation}
\end{remark}
Proof: Taking into account equation (\ref{eq4}) and properties of the covariation norm we obtain the following
\begin{eqnarray*}\|X_1(t)\|_{\alpha}&=&{\left\|\sum_{j=0}^{\infty} \left(a_1^{(j)} Z_1(t-j) + a_2^{(j)} Z_2(t-j)\right)\right\|}_{\alpha}\leq \sum_{j=0}^{\infty}{\left\|a_1^{(j)} Z_1(t-j) \right\|}_{\alpha}+
\sum_{j=0}^{\infty}{\left\|a_2^{(j)} Z_2(t-j)\right\|}_{\alpha}\nonumber\\
&=&\sum_{j=0}^{\infty}{\left(\left|a_1^{(j)}\right|\sigma_1+\left|a_2^{(j)}\right|\sigma_2\right)},
\end{eqnarray*}
where $\sigma_1$ and $\sigma_2$ are scale parameters for $\{Z_1(t)\}$ and $\{Z_2(t)\}$, respectively. One can easily observe $$\|X_1(t)\|_{\alpha}<\infty~~ \mbox{if} ~~\sum_{j=0}^{\infty}{\left(\left|a_1^{(j)}\right|+\left|a_2^{(j)}\right|\right)}<\infty.$$ 
Using the same reasoning for $X_2(t)$ given by Eq. (\ref{eq4}) we obtain the second condition and that asserts the result.
\begin{flushright}
$\Box$
\end{flushright}
In the further analysis we assume the model parameters satisfy condition given by (\ref{r2}).
\subsection{Cross-codifference}
Let us start our considerations from recalling the basic definition of the codifference for the symmetric $\alpha-$stable vector $(X,Y)$, which allows us to quantify the dependence between two random variables.

\begin{definition} \textup{\citep{Taqqu}} Let the random vector $(X,Y)$ be a bidimensional jointly and symmetric $\alpha-$stable i.e. its marginals have the following characteristic functions  $\phi_X(t)=\E \exp\{i t X\} =\exp\left\{-\sigma_X^\alpha |t|^\alpha\right\}$ and $\phi_Y(t)=\exp\left\{-\sigma_Y^\alpha |t|^\alpha\right\}$, respectively. Then the codifference between $X$ and $Y$ has the following form
$$
\mathrm{CD}(X,Y)=\sigma_{X-Y}^\alpha-\sigma_X^\alpha-\sigma_Y^\alpha.
\label{coddiference}
$$
\end{definition}
In literature \citep{Chechkin2015, nowicka} one can find also an alternative definition of codifference which uses the characteristic functions and can be defined for arbitrary random variables
$$
\mathrm{CD}(X,Y)=\log \E\exp\{i(X-Y)\}-\log \E\exp\{iX\} -\log \E\exp\{-iY\}.
\label{coddiference2}
$$
It is worth mentioning that both definitions of codifference in case of Gaussian random vector ($\alpha=2$) reduce to the usual covariance. Namely, we have that $\mathrm{CD}(X,Y)=\textrm{Cov}(X,Y)$ \citep{Taqqu,nowicka}. 

The codifference measure can be also used to quantify the time-dependence structure for stochastic processes. In \cite{Chechkin2015} an analogue of the auto-covariance was proposed in terms of the auto-codifference. 
\begin{definition}\textup{\citep{Chechkin2015}}  Let $\{X(t)\}$ be a stochastic process, then the auto-codifference is defined in the following way
$$
\mathrm{CD}(X(t),X(s))=\log \E\exp\{i(X(t)-X(s))\}-\log \E\exp\{iX(t)\}-\log \E\exp \{-iX(s)\}.
\label{autocoddiference}
$$
\end{definition}
For the stochastic processes that are stationary the auto-codifference depends only on the distance between arguments: $|t-s|$. Moreover, for an $\alpha-$stable process with $\alpha=2$ the auto-codifference reduces to the auto-covariance.\\

In the following definition we introduce the cross-codifference, new measure of dependence which can replace the cross-covariance defined only for second-order processes. The cross-codifference is an analogue of the classical cross-covariance. In probability and statistics, given two stochastic processes, the cross-covariance is a function that gives the covariance of one process with the other at pairs of time points. In our case, where we consider the bidimensional time series the cross-codifference will be calculated for components of bidimensional process, namely $X_1(t)$ and $X_2(t)$. 
\begin{definition}If $\{\mathbf{X}(t)\}=\{X_1(t),X_2(t)\}$  is a bidimensional process, then the cross-codifference is defined as
\begin{multline}
\mathrm{CD}(X_1(t),X_2(t+h))=\log \E \exp\{i(X_1(t)-X_2(t+h))\}-\log \E\exp\{iX_1(t)\} -\log \E \exp\{-iX_2(t+h)\}.
\label{eq5}
\end{multline}
\end{definition}
In the next Theorem we present the general formula for cross-codifference for general bidimensional VAR(1) model.
\begin{theorem}\label{th1}
For a general bidimensional VAR(1) process defined in (\ref{eq1}) and $h \geq 0$ the cross-codifference takes the following form
	\begin{enumerate}[(a)]
	\item 
	\begin{eqnarray}
	\label{eq6a}
&\mathrm{CD}\left(X_1\left(t\right),X_2\left(t+h\right)\right) &= \sum_{j=0}^{\infty}\log\left(\E\left(\exp\left\{i \left(\left(a_1^{\left(j\right)}-a_3^{\left(j+h\right)}\right) Z_1\left(t-j\right) + \left(a_2^{\left(j\right)}-a_4^{\left(j+h\right)}\right) Z_2\left(t-j\right)\right)\right\}\right)\right)\nonumber\\
& &-\sum_{j=0}^{\infty}\log\left(\E\left(\exp\left\{i\left(a_1^{\left(j\right)} Z_1\left(t-j\right) + a_2^{\left(j\right)} Z_2\left(t-j\right)\right)\right\}\right)\right)\\
& &-\sum_{j=0}^{\infty}\log\left(\E\left(\exp\left\{i\left(-a_3^{\left(j+h\right)} Z_1\left(t-j\right) - a_4^{\left(j+h\right)} Z_2\left(t-j\right)\right)\right\}\right)\right).\nonumber
\end{eqnarray}
	\item
	\begin{eqnarray}
	\label{eq6b}
&\mathrm{CD}(X_1(t),X_2(t-h))&=\sum_{j=h}^{\infty}\log\left(\E\left(\exp\left\{i \left(\left(a_1^{\left(j\right)}-a_3^{\left(j-h\right)}\right) Z_1\left(t-j\right) + \left(a_2^{\left(j\right)}-a_4^{\left(j-h\right)}\right) Z_2\left(t-j\right)\right)\right\}\right)\right)
\nonumber\\
& &-\sum_{j=h}^{\infty}\log\left(\E\left(\exp\left\{i\left(a_1^{\left(j\right)} Z_1\left(t-j\right) + a_2^{\left(j\right)} Z_2\left(t-j\right)\right)\right\}\right)\right)\\
& &-\sum_{j=h}^{\infty}\log\left(\E\left(\exp\left\{i\left(-a_3^{\left(j-h\right)} Z_1\left(t-j\right) - a_4^{\left(j-h\right)} Z_2\left(t-j\right)\right)\right\}\right)\right).\nonumber
		\end{eqnarray}
		\end{enumerate}
\end{theorem}
The proof of the above theorem is presented in Appendix A.

\section{Cross-codifference for bidimensional Gaussian VAR(1) model}
In this section we derive the analytical formula of cross-codifference function for the bidimensional Gaussian VAR(1) model with bidimensional noise defined in (\ref{eq7}). 
\begin{lemma}\label{l1} For a bidimensional Gaussian VAR(1) process with $\{\mathbf{Z}(t)\}$ given by (\ref{eq7}) and for $h \geq 0$ the cross-codifference has the following form
\begin{enumerate}[a)]
\item
\begin{equation}
\begin{split}
\mathrm{CD}\left(X_1\left(t\right),X_2\left(t+h\right)\right)=R_{11}\sum_{j=0}^{\infty}a_1^{\left(j\right)}a_3^{\left(j+h\right)}+R_{22}\sum_{j=0}^{\infty}a_2^{\left(j\right)}a_4^{\left(j+h\right)}+R_{12}\sum_{j=0}^{\infty}a_1^{\left(j\right)}a_4^{\left(j+h\right)}+R_{12}\sum_{j=0}^{\infty}a_2^{\left(j\right)}a_3^{\left(j+h\right)},
\end{split}
\label{eq9}
\end{equation}
\item
\begin{equation}
\begin{split}
\mathrm{CD}\left(X_1\left(t\right),X_2\left(t-h\right)\right)=R_{11}\sum_{j=0}^{\infty}a_1^{\left(j+h\right)}a_3^{\left(j\right)}+R_{22}\sum_{j=0}^{\infty}a_2^{\left(j+h\right)}a_4^{\left(j\right)}+R_{12}\sum_{j=0}^{\infty}a_1^{\left(j+h\right)}a_4^{\left(j\right)}+R_{12}\sum_{j=0}^{\infty}a_2^{\left(j+h\right)}a_3^{\left(j\right)}.
\end{split}
\label{eq9b}
\end{equation}
\end{enumerate}
\end{lemma}
Proof: {The proof follows directly from the general formulas for the cross-codifference of a bidimensional VAR(1) model given in (\ref{eq6a}) and (\ref{eq6b}) and from the formula for the characteristic function of $\{\mathbf{Z}(t)\}$ given in (\ref{eq8})}.
\begin{flushright}
$\Box$
\end{flushright}
\begin{remark}For bidimensional Gaussian VAR(1) models the cross-codifference is equal to the cross-covariance, namely
\[\mathrm{CD}(X_1(t),X_2(t+h))=\mathrm{Cov}(X_1(t),X_2(t+h)).\]
\end{remark}
Proof: If we take under consideration formula (\ref{eq4}), then one can easily calculate the cross-covariance, namely for $h \geq 0$ we have
\begin{eqnarray}
&\mathrm{Cov}(X_1(t),X_2(t+h))=\E\left(X_1(t)X_2(t+h)\right)=\nonumber\\&=\E\left(\sum_{j=0}^{\infty} (a_1^{(j)} Z_1(t-j) + a_2^{(j)} Z_2(t-j))\sum_{k=0}^{\infty} (a_3^{(k)} Z_1(t+h-k) + a_4^{(k)} Z_2(t+h-k))\right)=\nonumber\\
&= R_{11}\sum_{j=0}^{\infty} a_1^{(j)}a_3^{(j+h)}+R_{22}\sum_{j=0}^{\infty} a_2^{(j)}a_4^{(j+h)}+R_{12}\sum_{j=0}^{\infty} a_1^{(j)} a_4^{(j+h)}+R_{12}\sum_{j=0}^{\infty} a_2^{(j)} a_3^{(j+h)}=\nonumber\\
&=\mathrm{CD}(X_1(t),X_2(t+h)).
\end{eqnarray}
The similar reasoning implies that $\mathrm{Cov}(X_1(t),X_2(t-h))=\mathrm{CD}(X_1(t),X_2(t-h))$.
\begin{flushright}
$\Box$
\end{flushright}
\begin{example} As an example let us consider $R_{12}=0$, $a_2=a_3=0$ and $|a_1|<1$, $|a_4|<1$. Applying Lemma \ref{l1} we obtain
	\begin{equation*}
	\mathrm{CD}(X_1(t),X_2(t+h))=0
	\end{equation*}
and
	\begin{equation*}
	\mathrm{CD}(X_1(t),X_2(t-h))=0.
	\end{equation*}
	\end{example}
	\begin{example} As a second example let us consider $R_{12}\neq0$, $a_2=a_3=0$ and $|a_1|<1$, $|a_4|<1$. Applying Lemma \ref{l1} we obtain
	\begin{align*}
	\mathrm{CD}(X_1(t),X_2(t+h))=\sum_{j=0}^{\infty}a_1^{(j)}a_4^{(j+h)}R_{12}=R_{12}\sum_{j=0}^{\infty}a_1^{j}a_4^{j+h}=R_{12}a_4^{h}\sum_{j=0}^{\infty}{(a_1a_4)}^{j}=\frac{R_{12}a_4^{h}}{1-a_1a_4}
	\end{align*}
and
	\begin{align*}
	\mathrm{CD}(X_1(t),X_2(t-h))=\sum_{j=0}^{\infty}a_1^{(j+h)}a_4^{(j
		)}R_{12}=R_{12}\sum_{j=0}^{\infty}a_1^{j+h}a_4^{j}=R_{12}a_1^{h}\sum_{j=0}^{\infty}{(a_1a_4)}^{j}=\frac{R_{12}a_1^{h}}{1-a_1a_4}.
	\end{align*}
\end{example}

\section{Cross-codifference for bidimensional sub-Gaussian VAR(1) model}
In this section we derive the analytical formula of cross-codifference function for the sub-Gaussian bidimensional VAR(1) model with bidimensional noise $\{\mathbf{Z}(t)\}$ defined in (\ref{eq10}). 
\begin{lemma}\label{l2}
For a bidimensional sub-Gaussian VAR(1) model with $\{\mathbf{Z}(t)\}$ given by (\ref{eq10})  the cross-codifference  for $h \geq 0$ has the following form
\begin{enumerate}[a)]
\item 
\begin{eqnarray}
\begin{split}
&\mathrm{CD}\left(X_1\left(t\right),X_2\left(t+h\right)\right)={\left(\frac{1}{2}\right)}^{\alpha/2}\sum_{j=0}^{\infty}{\left|{\left(a_1^{\left(j\right)}\right)}^2 R_{11}+2a_1^{\left(j\right)}a_2^{\left(j\right)}R_{12}+{\left(a_2^{\left(j\right)}\right)}^2 R_{22}\right|}^{\alpha/2}\\
&+{\left(\frac{1}{2}\right)}^{\alpha/2}\sum_{j=0}^{\infty}{\left|{\left(a_3^{\left(j+h\right)}\right)}^2 R_{11}+2a_3^{\left(j+h\right)}a_4^{\left(j+h\right)}R_{12}+{\left(a_4^{\left(j+h\right)}\right)}^2R_{22}\right|}^{\alpha/2}\\
&-{\left(\frac{1}{2}\right)}^{\alpha/2}\sum_{j=0}^{\infty}{\left|{\left(a_1^{\left(j\right)}-a_3^{\left(j+h\right)}\right)}^2 R_{11}+2\left(a_1^{\left(j\right)}-a_3^{\left(j+h\right)}\right)\left(a_2^{\left(j\right)}-a_4^{\left(j+h\right)}\right)R_{12}+{\left(a_2^{\left(j\right)}-a_4^{\left(j+h\right)}\right)}^2R_{22}\right|}^{\alpha/2},
\end{split}
\label{eq14}
\end{eqnarray}
\item 
\begin{eqnarray}
\begin{split}
&\mathrm{CD}\left(X_1\left(t\right),X_2\left(t-h\right)\right)={\left(\frac{1}{2}\right)}^{\alpha/2}\sum_{j=0}^{\infty}{\left|{\left(a_1^{\left(j+h\right)}\right)}^2 R_{11}+2a_1^{\left(j+h\right)}a_2^{\left(j+h\right)}R_{12}+{\left(a_2^{\left(j+h\right)}\right)}^2R_{22}\right|}^{\alpha/2}\\
&+{\left(\frac{1}{2}\right)}^{\alpha/2}\sum_{j=0}^{\infty}{\left|{\left(a_3^{\left(j\right)}\right)}^2 R_{11}+2a_3^{\left(j\right)}a_4^{\left(j\right)}R_{12}+{\left(a_4^{\left(j\right)}\right)}^2R_{22}\right|}^{\alpha/2}\\
&-{\left(\frac{1}{2}\right)}^{\alpha/2}\sum_{j=0}^{\infty}{\left|{\left(a_1^{\left(j+h\right)}-a_3^{\left(j\right)}\right)}^2 R_{11}+2\left(a_1^{\left(j+h\right)}-a_3^{\left(j\right)}\right)\left(a_2^{\left(j+h\right)}-a_4^{\left(j\right)}\right)R_{12}+{\left(a_2^{\left(j+h\right)}-a_4^{\left(j\right)}\right)}^2R_{22}\right|}^{\alpha/2}.
\end{split}
\label{eq14b}
\end{eqnarray}
\end{enumerate}
\end{lemma}
{Proof: The proof follows directly from the general formulas for the cross-codifference of a bidimensional VAR(1) model given in (\ref{eq6a}) and (\ref{eq6b}) and from the formula for the characteristic function of $\{\mathbf{Z}(t)\}$ given in (\ref{eq11}).}	
\begin{flushright}
$\Box$
\end{flushright}

	\begin{example}\label{ex2}As an example let us consider $R_{12}\neq0$, $a_2=a_3=0$ and $|a_1|<1$, $|a_4|<1$. Applying Lemma \ref{l2} we obtain
	\begin{enumerate}[a)]
	\item
	\begin{equation*}
	\begin{split}
	&\mathrm{CD}\left(X_1\left(t\right),X_2\left(t+h\right)\right)=\\
	&={\left(\frac{1}{2}\right)}^{\alpha/2}\!\!\left(\sum_{j=0}^{\infty}{\left|{\left(a_1^{\left(j\right)}\right)}^2 R_{11}\right|}^{\alpha/2}\!\!\!+\!\sum_{j=0}^{\infty}{\left|{\left(a_4^{\left(j+h\right)}\right)}^2R_{22}\right|}^{\alpha/2}\!\!\!-\!\sum_{j=0}^{\infty}{\left|{\left(a_1^{\left(j\right)}\right)}^2 R_{11}\!-\!2 a_1^{\left(j\right)}a_4^{\left(j+h\right)}R_{12}\!+ \!{\left(a_4^{\left(j+h\right)}\right)}^2R_{22}\right|}^{\alpha/2}\right)\\
  &={\left(\frac{1}{2}\right)}^{\alpha/2}\!\!\left(\sum_{j=0}^{\infty}|a_1|^{\alpha j} {R_{11}}^{\alpha/2}+\sum_{j=0}^{\infty}|a_4|^{\alpha j+\alpha h}{R_{22}}^{\alpha/2}-\sum_{j=0}^{\infty}{\left|a_1^{2j} R_{11}-2 a_1^{j}a_4^{j+h}R_{12}+ a_4^{2j+2h}R_{22}\right|}^{\alpha/2}\right)\\
	&={\left(\frac{1}{2}\right)}^{\alpha/2}\!\!\left(\frac{R_{11}^{\alpha/2}}{1-|a_1|^\alpha} +\frac{R_{22}^{\alpha/2}|a_4|^{h\alpha}}{1-|a_4|^\alpha}-\sum_{j=0}^{\infty}{\left|a_1^{2j} R_{11}-2 a_1^{j}a_4^{j+h}R_{12}+a_4^{2j+2h}R_{22}\right|}^{\alpha/2}\right),
	\end{split}
	\end{equation*}
		\item
	\begin{equation*}
	\begin{split}
	&\mathrm{CD}\left(X_1\left(t\right),X_2\left(t-h\right)\right)=\\
	&={\left(\frac{1}{2}\right)}^{\alpha/2}\!\!\left(\sum_{j=0}^{\infty}{\left|{\left(a_1^{\left(j+h\right)}\right)}^2 R_{11}\right|}^{\alpha/2}\!\!\!+\!\sum_{j=0}^{\infty}{\left|{\left(a_4^{\left(j\right)}\right)}^2R_{22}\right|}^{\alpha/2}\!\!\!-\!\sum_{j=0}^{\infty}{\left|{\left(a_1^{\left(j+h\right)}\right)}^2 R_{11}\!-\!2 a_1^{\left(j+h\right)}a_4^{\left(j\right)}R_{12}\!+\!{\left(a_4^{\left(j\right)}\right)}^2R_{22}\right|}^{\alpha/2}\right)\\
	&={\left(\frac{1}{2}\right)}^{\alpha/2}\!\!\left(\sum_{j=0}^{\infty}|a_1|^{\alpha j+\alpha h} {R_{11}}^{\alpha/2}+\sum_{j=0}^{\infty}|a_4|^{\alpha j}{R_{22}}^{\alpha/2}-\sum_{j=0}^{\infty}{\left|a_1^{2j+2h} R_{11}-2 a_1^{j+h}a_4^{j}R_{12}+ a_4^{2j}R_{22}\right|}^{\alpha/2}\right)\\
	&={\left(\frac{1}{2}\right)}^{\alpha/2}\!\!\left(\frac{R_{11}^{\alpha/2}|a_1|^{h\alpha}}{1-|a_1|^\alpha} +\frac{R_{22}^{\alpha/2}}{1-|a_4|^\alpha}-\sum_{j=0}^{\infty}{\left|a_1^{2j+2h} R_{11}-2 a_1^{j+h}a_4^{j}R_{12}+a_4^{2j}R_{22}\right|}^{\alpha/2}\right).
	\end{split}
	\end{equation*}
	\end{enumerate}
\end{example}

\section{Simulations}
In this section we introduce an estimation method for cross-codifference from the experimental data. The idea is similar as presented in \cite{cod1} and is based on the replacement the theoretical characteristic function in the definition of the cross-codifference by its empirical equivalent. As the illustration we present exemplary trajectories of the considered time series together with a comparison of theoretical and empirical cross-codifference.\\

For a bidimensional time series $\{\textbf{X}(t)\}$ we define an estimator of cross-codifference as follows
\begin{multline}
\widehat{\mathrm{CD}}\left(X_1\left(t\right),X_2\left(s\right)\right)=\\=\log\left(\widehat{\phi}\left(1,-1,X_1\left(t\right),X_2\left(s\right)\right)\right)-\log\left(\widehat{\phi}\left(1,0,X_1\left(t\right),X_2\left(s\right)\right)\right)-\log\left(\widehat{\phi}\left(0,-1,X_1\left(t\right),X_2\left(s\right)\right)\right),
\end{multline}
where $\widehat{\phi}\left(u,v,X_1\left(t\right),X_2\left(s\right)\right)$ is an estimator of the following characteristic function
\begin{align*}
{\phi}\left(u,v,X_1\left(t\right),X_2\left(s\right)\right)=\E\left(\exp\left\{iu X_1\left(t\right)+i vX_2\left(s\right)\right\}\right).
\end{align*}
If we consider the realization of a stationary bidimensional time series  $\{\textbf{X}(t)\}$ denoted by $\{x^i_{k},k=1,2,\ldots,N,i=1,2\}$, the cross-codifference can be estimated from only one trajectory and in this case the estimator takes the following form
\begin{align}
\widehat{\mathrm{CD}}\left(X_1\left(t\right),X_2\left(t+k\right)\right)=\log\left(\widehat{\phi}\left(1,-1,k\right)\right)- \log\left(\widehat{\phi}\left(1,0,k\right)\right)-\log\left(\widehat{\phi}\left(0,-1,k\right)\right),
\label{ecodd}
\end{align}
where the empirical characteristic function $\widehat{\phi}(u,v,k)$ is given by:
\begin{equation}
\widehat{\phi}\left(u,v,k\right)= 
\begin{cases} \left(N-k\right)^{-1}\sum_{t=1}^{N-k}\exp\left(i\left(ux^1_{t}+vx^2_{t+k}\right)\right) & k\ge 0, \\
\left(N+k\right)^{-1}\sum_{t=1}^{N+k}\exp\left(i\left(ux^1_{t-k}+vx^2_{t}\right)\right) & k<0.
\end{cases}
\end{equation} 
The proposed estimator may be useful in quantifying the dependence structure for bidimensional processes. Similar as in one-dimensional case it can be also useful in the proper model recognition or estimation of appropriate parameters of the model \citep{Chechkin2015}.
\begin{figure}[h!]
	\centering
	\includegraphics[scale=.75]{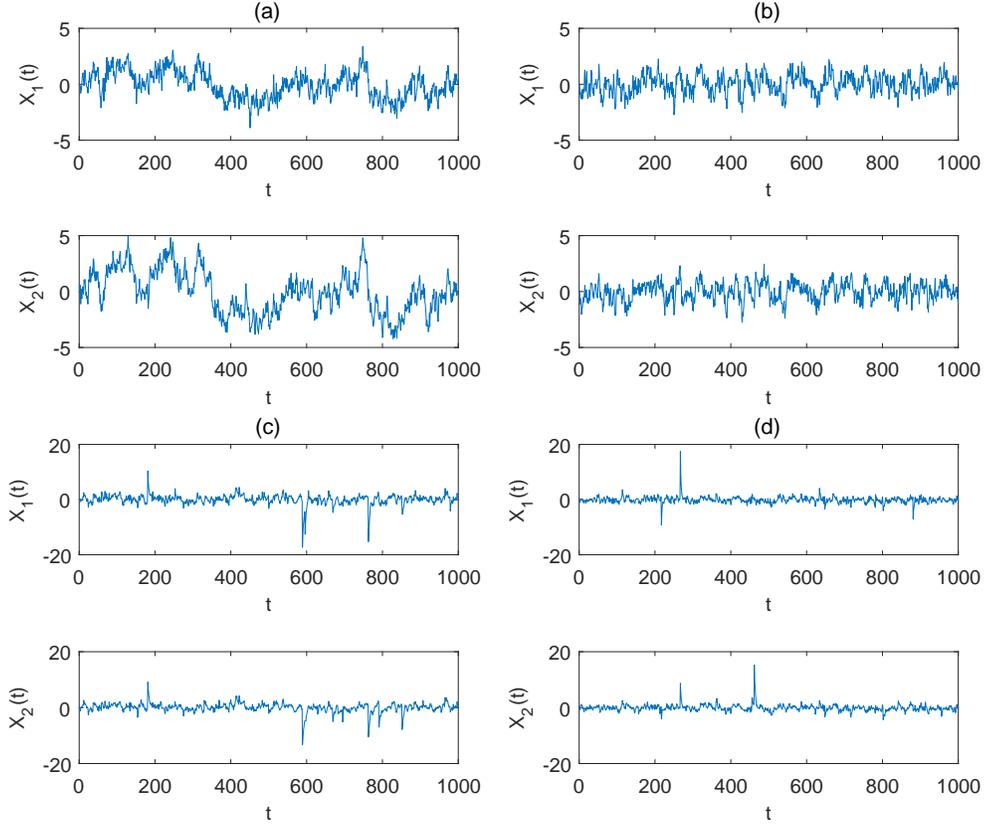}
	\caption{Exemplary trajectories of bidimensional time series: (a) bidimensional Gaussian VAR(1) model with $R_{11}=0.3$, $R_{22}=0.3$, $R_{12}=0.2$, $a_1=0.6$, $a_2=0.2$, $a_3=0.1$, $a_4=0.9$, (b) bidimensional Gaussian VAR(1) model with $R_{11}=0.3$, $R_{22}=0.3$, $R_{12}=0.2$, $a_1=0.6$, $a_2=0$, $a_3=0$, $a_4=0.9$, (c) bidimensional sub-Gaussian VAR(1) model with $\alpha=1.7$, $R_{11}=0.4$, $R_{22}=0.3$, $R_{12}=0.3$, $a_1=0.6$, $a_2=0.2$, $a_3=0.1$, $a_4=0.7$, (d) bidimensional sub-Gaussian VAR(1) model with $\alpha=1.7$, $R_{11}=0.4$, $R_{22}=0.3$, $R_{12}=0.3$, $a_1=0.6$, $a_2=0$, $a_3=0$, $a_4=0.7$.}
	\label{fig1}
\end{figure}
\begin{figure}[h!]
	\centering
	\includegraphics[scale=.5]{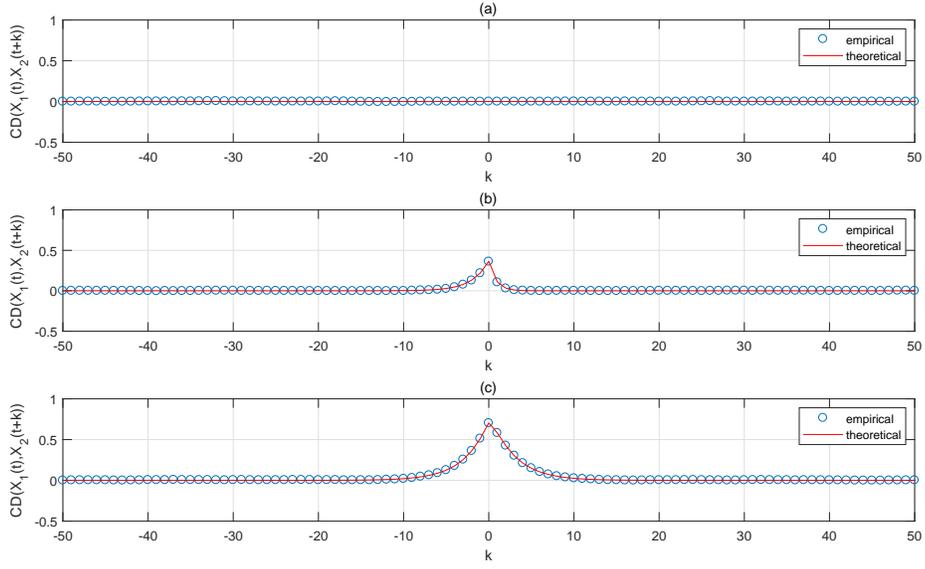}
	\caption{Estimator of cross-codifference for the bidimensional Gaussian VAR(1) model and the theoretical values given by (\ref{eq9}) and (\ref{eq9b}): (a) $R_{11}=0.5$, $R_{22}=0.5$, $R_{12}=0$, $a_1=0.6$, $a_2=0$, $a_3=0$, $a_4=0.3$, (b) $R_{11}=0.5$, $R_{22}=0.5$, $R_{12}=0.3$, $a_1=0.6$, $a_2=0$, $a_3=0$, $a_4=0.3$, (c) $R_{11}=0.5$, $R_{22}=0.5$, $R_{12}=0.3$, $a_1=0.6$, $a_2=0.1$, $a_3=0.4$, $a_4=0.3$.}
	\label{fig2}
\end{figure}
\begin{figure}[h!]
	\centering
	\includegraphics[scale=.5]{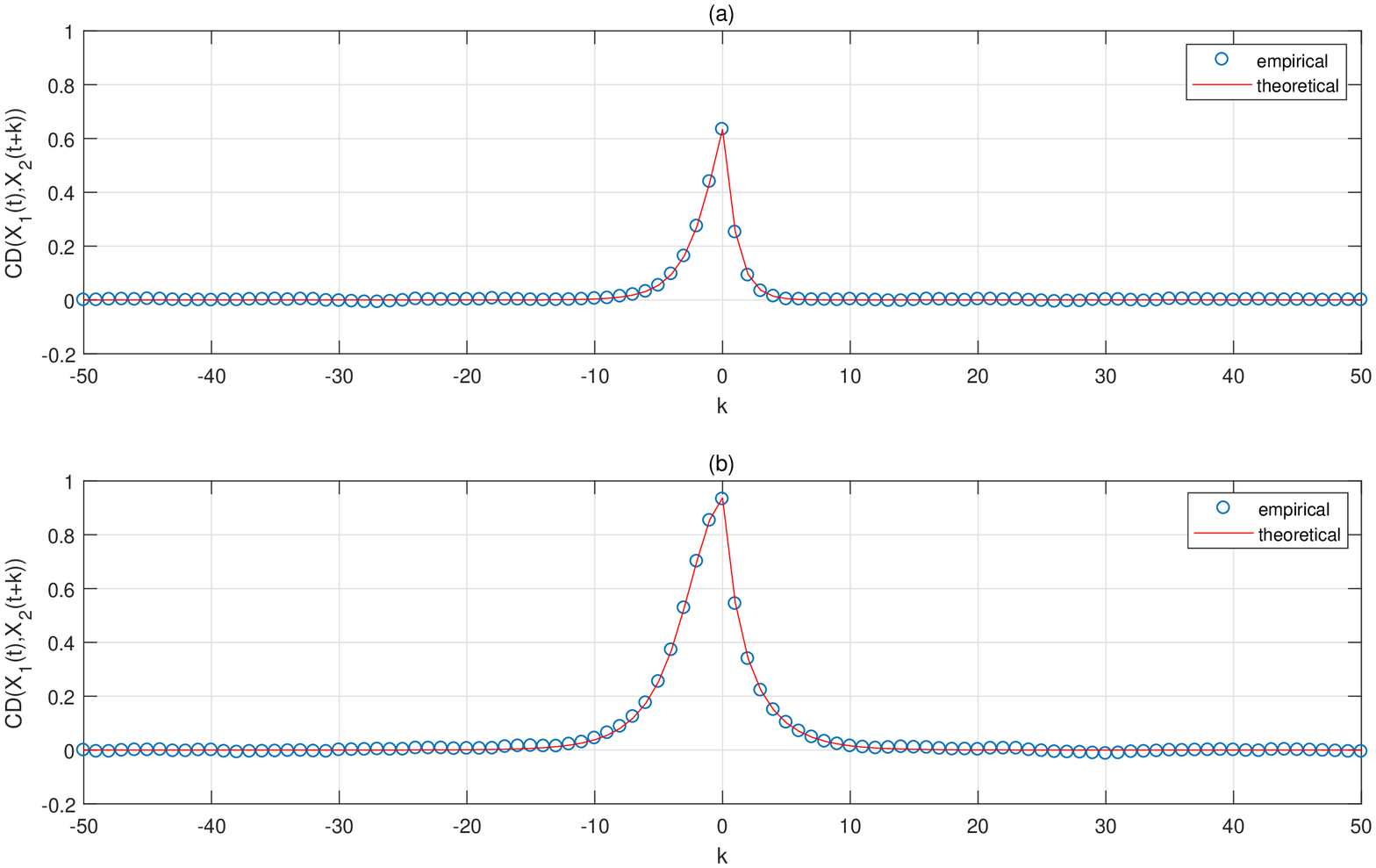}
	\caption{Estimator of cross-codifference for the bidimensional sub-Gaussian VAR(1) model and the theoretical values given by (\ref{eq14}) and (\ref{eq14b}): (a) $\alpha=1.5$, $R_{11}=0.4$, $R_{22}=0.3$, $R_{12}=0.3$, $a_1=0.6$, $a_2=0$, $a_3=0$, $a_4=0.4$, (b) $\alpha=1.5$, $R_{11}=0.4$, $R_{22}=0.3$, $R_{12}=0.3$, $a_1=0.6$, $a_2=0.3$, $a_3=0.1$, $a_4=0.4$.}
	\label{fig3}
\end{figure}

In Figure \ref{fig1} we show sample paths of the considered time series. On the top panels (a) and (b) of Figure \ref{fig1} we present the exemplary trajectories of bidimensional Gaussian VAR(1) time series, whereas on the bottom panels (c) and (d) we present the exemplary trajectories of bidimensional sub-Gaussian VAR(1) time series. For the model with infinite variance one can see the difference in the amplitude of the observations. Now, the next step is to verify the theoretical formulas for the cross-codifference given in the Section 3 and in the Section 4. In order to perform the comparison, we generate sample trajectories of the considered VAR(1) models. Using simulated data we calculate the empirical cross-codifference given in (\ref{ecodd}) and we plot it together with the corresponding theoretical formulas. The results for the Gaussian model are presented in Figure \ref{fig2}: (a) corresponds to Example 3.1 presented in Section 3, (b) corresponds to Example 3.2 presented in Section 3 and (c) is a general example. The results for the sub-Gaussian VAR(1) time series are presented in Figure \ref{fig3}: (a) corresponds to Example 4.1 presented in Section 4 and (b) is a general example. To calculate the theoretical values we use the formulas given in (\ref{eq9}), (\ref{eq9b}), (\ref{eq14}), (\ref{eq14b}) by taking a sum over $j$ from $0$ to $50$. In all cases one observes almost perfect agreement between the empirical and theoretical results.

\section{Estimation}

In this section we demonstrate how the theoretical results presented in the previous parts of the paper can be applied to the estimation of the considered model parameters. As the example, we introduce an estimation procedure for the parameters of bidimensional sub-Gaussian VAR(1) model. The method is based on the formula for the cross-codifference presented in Section 4, Example \ref{ex2}. We discuss here the case when the components of bidimensional time series are dependent only in the sense of the bidimensional noise. This is the case when. $a_2 = a_3 = 0$. More precisely, we consider Example \ref{ex2} with $0<a_1<1$ and $0<a_4<1$.
\begin{figure}[h!]
	\centering
	\includegraphics[scale=.6]{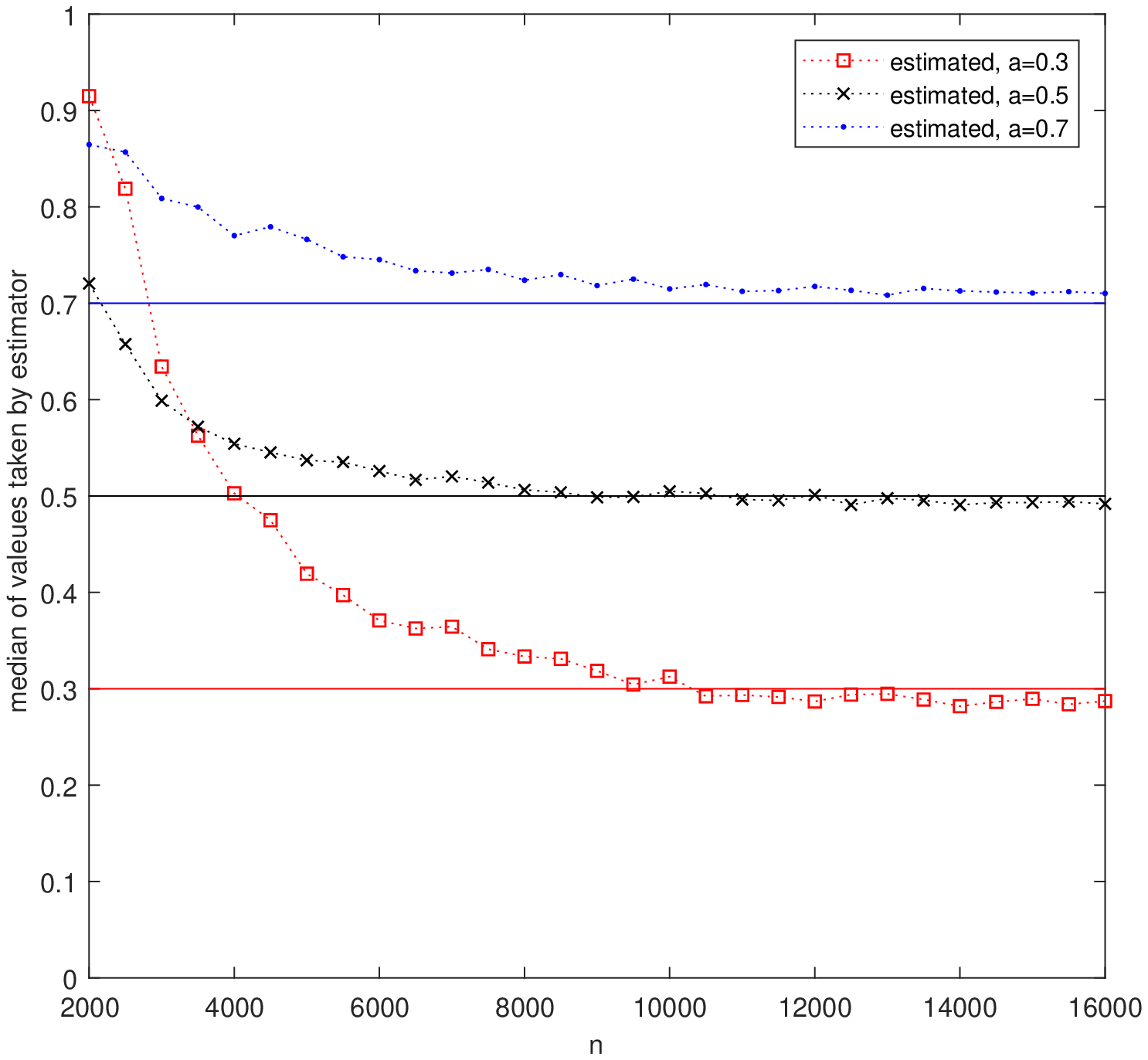}
	\caption{Medians of estimator $\widehat{a}$ for trajectories of various lengths $n$. In order to obtain the median of the estimators we simulated $1000$ trajectories for each length. }
	\label{fig4}
\end{figure}

At first, let us take $a_1=a_2=a$. Moreover, we assume $\alpha>1$. In this case, according to Example \ref{ex2} we have:
\[CD(X_1(t),X_2(t+h))=\frac{{\left(\frac{1}{2}\right)}^{\alpha/2}}{1-|a|^\alpha}\!\!\left(R_{11}^{\alpha/2} +R_{22}^{\alpha/2}|a|^{h\alpha}-\left|R_{11}-2 a^{h}R_{12}+a^{2h}R_{22}\right|^{\alpha/2}\right)\]
\[CD(X_1(t),X_2(t-h))=\frac{{\left(\frac{1}{2}\right)}^{\alpha/2}}{1-|a|^\alpha}\!\!\left({R_{11}^{\alpha/2}|a|^{h\alpha}} +{R_{22}^{\alpha/2}}-{\left|a^{2h} R_{11}-2 a^{h}R_{12}+R_{22}\right|}^{\alpha/2}\right).\]
One can show by using using L'H\^opital rule that the following holds for $\alpha>1$: 
\[\lim_{x\rightarrow 0}\frac{A^{{\alpha/2}}|x|^{\alpha} +B^{{\alpha/2}}-{\left|Ax^{2}-2 Cx+B\right|}^{\alpha/2}}{|x| }=D,\]
where $A,B,C,D$ are some constants. Therefore in the considered case, it can be proven that the cross-codifference behaves asymptotically as:
 $$CD(X_1(t),X_2(t+h)) \sim c_1a^h~ \mbox{and}~~ CD(X_1(t),X_2(t-h)) \sim c_2a^h,$$ where $c_1$ and $c_2$ are some constants.

In our methodology for bidimensional vector of observations we calculate the empirical cross-codifference $\widehat{\rm{CD}}(X_1(t),X_2(t+h))$ for $h>0$ and by comparing its values with the asymptotic formula we estimate the unknown parameters $a$ and $c_1$ using least squares method:
\begin{equation}
(\widehat{c_1},\widehat{a}) = \min\limits_{c_1,a}\sum_{i}(c_1a^{h_i}-\widehat{\rm{CD}}(X_1(t),X_2(t+h_i)))^2.
\end{equation}
To verify the effectiveness of the estimator we perform the Monte Carlo study. We generate 1000 trajectories of considered bidimensional time series and for each trajectory we estimate $\widehat{a}$. Then, we calculate the median of values taken by estimator. We repeat the simulations for the trajectories of various lengths. Exemplary results are presented in Figure \ref{fig4} where the convergence of the estimator is visible. The larger the length of a trajectory, the closer to the theoretical value is the outcome. 

As the second example, let us consider $a_1 \neq a_4$. In this case we can also observe the behaviour similar to the previous instance. Using the same methodology, one can show the cross-codifference converges to zero in the following manner:
$$CD(X_1(t),X_2(t+h)) \sim d_1a_4^h~~ \mbox{and}~~ CD(X_1(t),X_2(t-h)) \sim d_2a_1^h,$$ where $d_1$ and $d_2$ are some constants.
\begin{figure}[h!]
	\centering
	\includegraphics[scale=.65]{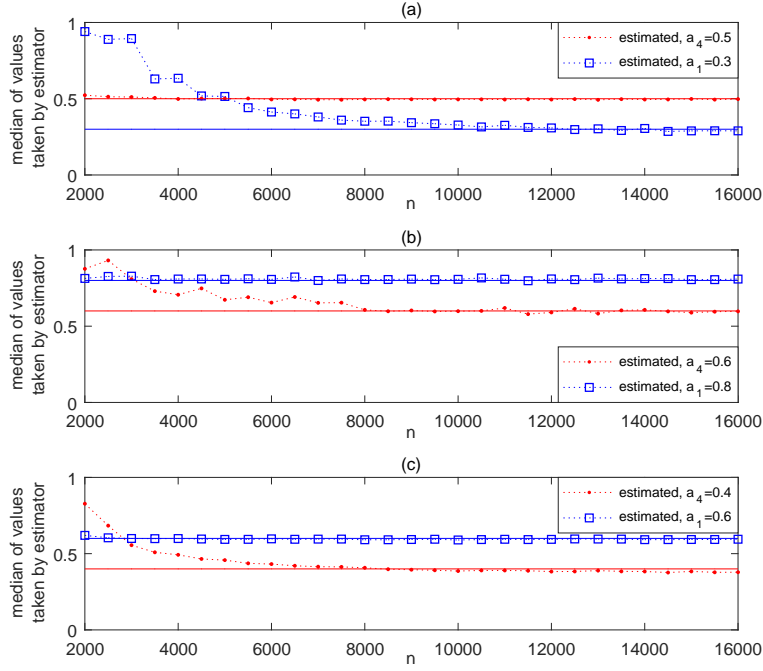}
	\caption{Medians of estimators $\widehat{a_1}$ and $\widehat{a_4}$ for trajectories of various lengths $n$. In order to obtain the median of the estimators we simulated $1000$ trajectories for each length.}
	\label{fig5}
\end{figure}
Similarly to the first case, in order to estimate the unknown parameters  $a_1$ and $a_4$ for bidimensional vector of observations we calculate the empirical cross-codifference $\widehat{\rm{CD}}(X_1(t),X_2(t+h))$ and $\widehat{\rm{CD}}(X_1(t),X_2(t-h))$ for $h>0$ and we compare the values with the asymptotic formulas using the least squares method:
\begin{equation}
(\widehat{d_2},\widehat{a_4}) = \min\limits_{d_2,a_4}\sum_{i}(d_2a_4^{h_i}-\widehat{\rm{CD}}(X_1(t),X_2(t+h_i)))^2
\end{equation}
and
\begin{equation}
(\widehat{d_1},\widehat{a_1}) = \min\limits_{d_1,a_1}\sum_{i}(d_1a_1^{h_i}-\widehat{\rm{CD}}(X_1(t),X_2(t-h_i)))^2.
\end{equation}
The effectiveness of the estimator is verified using Monte Carlo simulations and the exemplary results are presented in Figure \ref{fig5}. One can observe that for the parameter taking larger value the estimator converges to its theoretical equivalent much faster than for the other one.

\begin{figure}[h!]
	\centering
	\includegraphics[scale=.65]{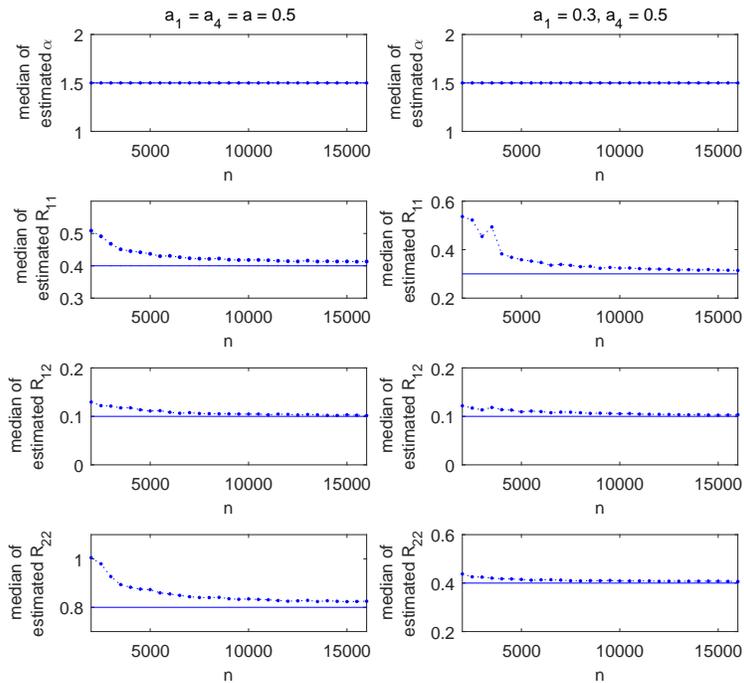}
	\caption{Medians of estimators $\widehat{\alpha}$, $\widehat{R}_{11}$, $\widehat{R}_{12}$ and $\widehat{R}_{22}$ for trajectories of various lengths $n$. Left panels correspond to the case when $a_1=a_4=a$ with $a=0.5$, $\alpha=1.5$, $R_{11}=0.4$, $R_{12}=0.1$ and $R_{22}=0.8$ and right panels correspond to the case when $a_1 \neq a_4$ with $a_1=0.3$, $a_4=0.5$, $\alpha=1.5$, $R_{11}=0.3$, $R_{12}=0.1$ and $R_{22}=0.4$. In order to obtain the median of the estimators we simulated $1000$ trajectories for each length.}
	\label{fig6}
\end{figure}

In order to estimate the remaining parameters related to the noise, after estimation of the VAR(1) model parameters we propose to extract the noise from the data by applying the inverse filter to each component of bidimensional time series by using estimated values of VAR(1) model. In the next step, we use the methodology introduced in \cite{sub4} where the estimation procedure based on the distance between the empirical and theoretical characteristic function of bivariate sub-Gaussian vectors is presented. The exemplary results obtained via Monte Carlo simulations are presented in Figure \ref{fig6}. The significant impact on the values taken by these estimators has the goodness of $a_1$ and $a_4$ estimation.
\section{Conclusions}
In this paper we propose a cross-codifference measure as a tool to identify the dependence structure between spatial components of multivariate stochastic processes with non-Gaussian innovations. First, we have established the general form of the solutions of bidimensional VAR(1) time series with Gaussian and sub-Gaussian innovations. The relation between the cross-codifference and cross-covariance in case of VAR(1) model with bidimensional Gaussian innovations have been obtained. The main practical result of this work is a derivation of the cross-codifference for the bidimensional VAR(1) model with sub-Gaussian noise and a proposition of its estimation technique. Moreover, we have proposed a new estimation technique for VAR(1) model parameters using the cross-codifference. It was also shown based on the simulated data that the introduced cross-dependence measure can be a useful tool in modelling of the data that are characterized by the dependence between spatial components in time. The study carried out in this paper open up a new areas of interest in context of the cross-dependence of multivariate models and estimation of their parameters. 

\section*{Acknowledgments}
AG and AW would like to acknowledge a support of National Center of Science Opus Grant No. 2016/21/B/ST1/00929 "Anomalous diffusion processes and their applications in real data modelling".
\section*{Appendix A}
Proof of Theorem \ref{th1}.\\
\leavevmode 
\begin{enumerate}[(a)]
\item First we compute each component of (\ref{eq5}) separately using the formulas for $X_1(t)$ and $X_2(t+h)$ given in (\ref{eq4}). We assume $h>0$.
\begin{align}
&\log\left(\E\left(\exp\left\{iX_1(t)\right\}\right)\right) 
=\log\biggl(\mathrm{E}\biggl(\exp\biggl\{i\sum_{j=0}^{\infty} \left(a_1^{(j)} Z_1(t-j) + a_2^{(j)} Z_2(t-j)\right)\biggr\}\biggr)\biggr)=\nonumber\\
&=\log\left(\E\left(\prod_{j=0}^{\infty}\exp\left\{i(a_1^{(j)} Z_1(t-j) + a_2^{(j)} Z_2(t-j))\right\}\right)\right)=\nonumber\\
&=\log\bigl(\prod_{j=0}^{\infty}\E(\exp\{i(a_1^{(j)} Z_1(t-j) + a_2^{(j)} Z_2(t-j))\})\bigr)=\nonumber\\
&=\sum_{j=0}^{\infty}\log\bigl(\E(\exp\{i(a_1^{(j)} Z_1(t-j) + a_2^{(j)} Z_2(t-j))\})\bigr)
\label{proof1}
\end{align}
\begin{align}
&\log\bigr(E(\exp\{-iX_2(t+h)\})\bigl)=\log\bigr(E(\exp\{-i\sum_{k=0}^{\infty} (a_3^{(k)} Z_1(t+h-k) + a_4^{(k)} Z_2(t+h-k))\})\bigl)= \nonumber\\
&=\log\bigr(E(\exp\{-i\sum_{j=-h}^{\infty} (a_3^{(j+h)} Z_1(t-j) + a_4^{(j+h)} Z_2(t-j))\})\bigl)= \nonumber\\
&=\log\bigr(E(\prod_{j=-h}^{\infty}\exp\{-i(a_3^{(j+h)} Z_1(t-j) + a_4^{(j+h)} Z_2(t-j))\})\bigl)= \nonumber\\
&=\log\bigr(\prod_{j=-h}^{\infty}E(\exp\{-i(a_3^{(j+h)} Z_1(t-j) + a_4^{(j+h)} Z_2(t-j))\})\bigl)= \nonumber\\
&=\sum_{j=-h}^{\infty}\log\bigr(E(\exp\{-i(a_3^{(j+h)} Z_1(t-j) + a_4^{(j+h)} Z_2(t-j))\})\bigl)= \nonumber\\
&=\sum_{j=-h}^{-1}\log\bigr(E(\exp\{i(-a_3^{(j+h)} Z_1(t-j) - a_4^{(j+h)} Z_2(t-j))\})\bigl)+ \nonumber\\
&+\sum_{j=0}^{\infty}\log\bigr(E(\exp\{i(-a_3^{(j+h)} Z_1(t-j) - a_4^{(j+h)} Z_2(t-j))\})\bigl)
\label{proof2}
\end{align}
Let us consider the random variable $X_1(t)-X_2(t+h)$.
\begin{align*}
&X_1(t)-X_2(t+h)= \sum_{j=0}^{\infty} (a_1^{(j)} Z_1(t-j) + a_2^{(j)} Z_2(t-j)) - \sum_{k=0}^{\infty} (a_3^{(k)} Z_1(t+h-k) + a_4^{(k)} Z_2(t+h-k)) =\\
&=\sum_{j=0}^{\infty} (a_1^{(j)} Z_1(t-j) + a_2^{(j)} Z_2(t-j)) - \sum_{j=-h}^{\infty} (a_3^{(j+h)} Z_1(t-j) + a_4^{(j+h)} Z_2(t-j)) =\\
&=\sum_{j=0}^{\infty} (a_1^{(j)} Z_1(t-j) + a_2^{(j)} Z_2(t-j)) - \sum_{j=-h}^{-1} (a_3^{(j+h)} Z_1(t-j) + a_4^{(j+h)} Z_2(t-j)) -\\
&- \sum_{j=0}^{\infty} (a_3^{(j+h)} Z_1(t-j) + a_4^{(j+h)} Z_2(t-j)) =\\
&= \sum_{j=0}^{\infty} \underbrace{((a_1^{(j)}-a_3^{(j+h)}) Z_1(t-j) + (a_2^{(j)}-a_4^{(j+h)}) Z_2(t-j))}_{A^{*}_{j}} - \sum_{j=-h}^{-1} \underbrace{(a_3^{(j+h)} Z_1(t-j) + a_4^{(j+h)} Z_2(t-j))}_{B^{*}_{j}}
\end{align*}
\begin{align}
&\log\bigr(E(\exp\{i(X_1(t)-X_2(t+h))\})\bigl)=\log\bigr(E(\exp\{i(\sum_{j=0}^{\infty} A^{*}_j - \sum_{j=-h}^{-1} B^{*}_j)\})\bigl)= \nonumber\\
&=\log\bigr(E(\exp\{i\sum_{j=0}^{\infty} A^{*}_j\}\ \exp\{-i \sum_{j=-h}^{-1} B^{*}_j)\})\bigl)=\log\bigr(E(\exp\{i\sum_{j=0}^{\infty} A^{*}_j\})\ E(\exp\{-i \sum_{j=-h}^{-1} B^{*}_j\}\bigl)= \nonumber\\
&=\log\bigr(E(\prod_{j=0}^{\infty}\exp\{i A^{*}_j\})\ E(\prod_{j=-h}^{-1}\exp\{-i B^{*}_j\}\bigl)= \nonumber\\
&=\log\bigr(E(\prod_{j=0}^{\infty}\exp\{i ((a_1^{(j)}-a_3^{(j+h)}) Z_1(t-j) + (a_2^{(j)}-a_4^{(j+h)}) Z_2(t-j))\})\nonumber\\
& E(\prod_{j=-h}^{-1}\exp\{-i (a_3^{(j+h)} Z_1(t-j) + a_4^{(j+h)} Z_2(t-j))\}\bigl)=\nonumber\\
&=\log\bigr(\prod_{j=0}^{\infty}E(\exp\{i ((a_1^{(j)}-a_3^{(j+h)}) Z_1(t-j) + (a_2^{(j)}-a_4^{(j+h)}) Z_2(t-j))\})\nonumber\\ 
&\prod_{j=-h}^{-1}E(\exp\{-i (a_3^{(j+h)} Z_1(t-j) + a_4^{(j+h)} Z_2(t-j))\}\bigl)=\nonumber\\
&=\log\bigr(\prod_{j=0}^{\infty}E(\exp\{i ((a_1^{(j)}-a_3^{(j+h)}) Z_1(t-j) + (a_2^{(j)}-a_4^{(j+h)}) Z_2(t-j))\})\bigr)+\nonumber\\
&+\log \bigl( \prod_{j=-h}^{-1}E(\exp\{-i (a_3^{(j+h)} Z_1(t-j) + a_4^{(j+h)} Z_2(t-j))\}\bigl)=\nonumber\\
&=\sum_{j=0}^{\infty}\log\bigr(E(\exp\{i ((a_1^{(j)}-a_3^{(j+h)}) Z_1(t-j) + (a_2^{(j)}-a_4^{(j+h)}) Z_2(t-j))\})\bigr)+\nonumber\\
&+\sum_{j=-h}^{-1}\log \bigl( E(\exp\{i (-a_3^{(j+h)} Z_1(t-j) - a_4^{(j+h)} Z_2(t-j))\}\bigl)
\label{proof3}
\end{align}
Finally, taking into account (\ref{proof1}), (\ref{proof2}) and (\ref{proof3}), we get:
\begin{align*}
&\mathrm{CD}(X_1(t),X_2(t+h))= -\sum_{j=0}^{\infty}\log\bigr(E(\exp\{i(a_1^{(j)} Z_1(t-j) + a_2^{(j)} Z_2(t-j))\})\bigl)\nonumber\\
&-\sum_{j=0}^{\infty}\log\bigr(E(\exp\{i(-a_3^{(j+h)} Z_1(t-j) - a_4^{(j+h)} Z_2(t-j))\})\bigl)\\
&+\sum_{j=0}^{\infty}\log\bigr(E(\exp\{i ((a_1^{(j)}-a_3^{(j+h)}) Z_1(t-j) + (a_2^{(j)}-a_4^{(j+h)}) Z_2(t-j))\})\bigr).\nonumber
\end{align*}
\item Again, we compute each component of (\ref{eq5}) separately using the formulas for $X_1(t)$ and $X_2(t-h)$ given in (\ref{eq4}). We assume $h>0$. The first component is given by (\ref{proof1}).
\begin{align}
&\log\bigr(E(\exp\{-iX_2(t-h)\})\bigl)=\log\bigr(E(\exp\{-i\sum_{k=0}^{\infty} (a_3^{(k)} Z_1(t-h-k) + a_4^{(k)} Z_2(t-h-k))\})\bigl)=\nonumber\\
&=\log\bigr(E(\exp\{-i\sum_{j=h}^{\infty} (a_3^{(j-h)} Z_1(t-j) + a_4^{(j-h)} Z_2(t-j))\})\bigl)=\nonumber\\
&=\log\bigr(E(\prod_{j=h}^{\infty}\exp\{-i(a_3^{(j-h)} Z_1(t-j) + a_4^{(j-h)} Z_2(t-j))\})\bigl)=\nonumber\\
&=\log\bigr(\prod_{j=h}^{\infty}E(\exp\{-i(a_3^{(j-h)} Z_1(t-j) + a_4^{(j-h)} Z_2(t-j))\})\bigl)=\nonumber\\
&=\sum_{j=h}^{\infty}\log\bigr(E(\exp\{-i(a_3^{(j-h)} Z_1(t-j) + a_4^{(j-h)} Z_2(t-j))\})\bigl)
\label{proof4}
\end{align}
Let us consider the random variable $X_1(t)-X_2(t-h)$.
\begin{align*}
X_1(t)-X_2(t-h)= \sum_{j=0}^{\infty} (a_1^{(j)} Z_1(t-j) + a_2^{(j)} Z_2(t-j)) - \sum_{k=0}^{\infty} (a_3^{(k)} Z_1(t-h-k) + a_4^{(k)} Z_2(t-h-k)) =\\
=\sum_{j=0}^{\infty} (a_1^{(j)} Z_1(t-j) + a_2^{(j)} Z_2(t-j)) - \sum_{j=h}^{\infty} (a_3^{(j-h)} Z_1(t-j) + a_4^{(j-h)} Z_2(t-j)) =\\
=\sum_{j=0}^{h-1} \underbrace{(a_1^{(j)} Z_1(t-j) + a_2^{(j)} Z_2(t-j))}_{C_{j}^{*}} + \sum_{j=h}^{\infty} \underbrace{((a_1^{(j)}-a_3^{(j-h)}) Z_1(t-j) + (a_2^{(j)}-a_4^{(j-h)}) Z_2(t-j))}_{D_{j}^{*}}
\end{align*}
\begin{align}
&\log\bigr(E(\exp\{i(X_1(t)-X_2(t-h))\})\bigl)=\log\bigr(E(\exp\{i(\sum_{j=0}^{h-1} C^{*}_j + \sum_{j=h}^{\infty} D^{*}_j)\})\bigl)=\nonumber\\
&=\log\bigr(E(\exp\{i\sum_{j=0}^{h-1} C^{*}_j\}\ \exp\{i \sum_{j=h}^{\infty} D^{*}_j)\})\bigl)=\log\bigr(E(\exp\{i\sum_{j=0}^{h-1} C^{*}_j\})\ E(\exp\{-i \sum_{j=h}^{\infty} D^{*}_j\}\bigl)=\nonumber\\
&=\log\bigr(E(\prod_{j=0}^{h-1}\exp\{i C^{*}_j\})\ E(\prod_{j=h}^{\infty}\exp\{-i D^{*}_j\}\bigl)=\nonumber\\
&=\log\bigr(E(\prod_{j=0}^{h-1}\exp\{i ((a_1^{(j)} Z_1(t-j) + a_2^{(j)} Z_2(t-j)))\})\nonumber\\
& E(\prod_{j=h}^{\infty}\exp\{i ((a_1^{(j)}-a_3^{(j-h)}) Z_1(t-j) + (a_2^{(j)}-a_4^{(j-h)}) Z_2(t-j))\})\bigl)=\nonumber\\
&=\log\bigr(\prod_{j=0}^{h-1}E(\exp\{i ((a_1^{(j)} Z_1(t-j) + a_2^{(j)} Z_2(t-j)))\})\nonumber\\
& \prod_{j=h}^{\infty}E(\exp\{i ((a_1^{(j)}-a_3^{(j-h)}) Z_1(t-j) + (a_2^{(j)}-a_4^{(j-h)}) Z_2(t-j))\})\bigl)=\nonumber\\
&=\log\bigr(\prod_{j=0}^{h-1}E(\exp\{i ((a_1^{(j)} Z_1(t-j) + a_2^{(j)} Z_2(t-j)))\})\bigl)+\nonumber\\ 
&\log\bigr(\prod_{j=h}^{\infty}E(\exp\{i ((a_1^{(j)}-a_3^{(j-h)}) Z_1(t-j) + (a_2^{(j)}-a_4^{(j-h)}) Z_2(t-j))\})\bigl)=\nonumber\\
&=\sum_{j=0}^{h-1}\log\bigr(E(\exp\{i ((a_1^{(j)} Z_1(t-j) + a_2^{(j)} Z_2(t-j)))\})\bigl)+\nonumber\\ 
&\sum_{j=h}^{\infty}\log\bigr(E(\exp\{i ((a_1^{(j)}-a_3^{(j-h)}) Z_1(t-j) + (a_2^{(j)}-a_4^{(j-h)}) Z_2(t-j))\})\bigl)
\label{proof5}
\end{align}
Finally, taking into account (\ref{proof1}), (\ref{proof4}) and (\ref{proof5}), we get:
\begin{align*}
&\mathrm{CD}(X_1(t),X_2(t-h))= -\sum_{j=h}^{\infty}\log\bigr(E(\exp\{i(a_1^{(j)} Z_1(t-j) + a_2^{(j)} Z_2(t-j))\})\bigl)\nonumber\\
&-\sum_{j=h}^{\infty}\log\bigr(E(\exp\{i(-a_3^{(j-h)} Z_1(t-j) - a_4^{(j-h)} Z_2(t-j))\})\bigl)\\
&+\sum_{j=h}^{\infty}\log\bigr(E(\exp\{i ((a_1^{(j)}-a_3^{(j-h)}) Z_1(t-j) + (a_2^{(j)}-a_4^{(j-h)}) Z_2(t-j))\})\bigr).\nonumber
\end{align*}
\begin{flushright}
	$\Box$
\end{flushright}
\end{enumerate}

\bibliography{moja_bib}

\end{document}